\documentclass[reqno]{amsart}
\usepackage{amssymb}

\usepackage{amssymb,epsfig,graphics}
\usepackage{amsmath}
\usepackage{amscd}
\usepackage{verbatim}

\input xypic

\newtheorem{thm}{Theorem}[section]
\newtheorem{prop}[thm]{Proposition}

\newtheorem{coro}[thm]{Corollary}

\theoremstyle{definition}

\newtheorem{example}{Example}

\def\R{\mathbb R}

\def\r{\rangle}
\def\l{\langle}

\def\a{\alpha}

\def\l{\langle}
\def\r{\rangle}

\begin{document}

\title[]{Faces of Platonic solids in all dimensions}

\author{M.~Szajewska}

\maketitle

\noindent
Centre de recherches math\'ematiques, Universit\'e de Montr\'eal, C.~P.~6128 -- Centre ville, Montr\'eal, H3C\,3J7, Qu\'ebec, Canada.

\noindent
Permanent Address: Institute of Mathematics, University of Bialystok, Akademicka~2, PL-15-267, Bialystok, Poland.

\vspace{10pt}

\noindent
E-mail: m.szajewska@math.uwb.edu.pl

\begin{abstract}\

This paper considers Platonic solids/polytopes in the real Euclidean space $\R^n$ of dimension $3\leq n<\infty$. The Platonic solids/polytopes are described together with their faces of dimensions $0\leq d\leq n-1$. Dual pairs of Platonic polytopes are considered in parallel. The underlying finite Coxeter groups are those of simple Lie algebras of types $A_n$, $B_n$, $C_n$, $F_4$ and of non-crystallographic Coxeter groups $H_3$, $H_4$.

Our method consists in recursively decorating the appropriate Coxeter-Dynkin diagram. Each recursion step provides the essential information about faces of a specific dimension. If, at each recursion step, all of the faces are in the same Coxeter group orbit, i.e. are identical, the solid is called Platonic.

\end{abstract}
\smallskip



\section{Introduction}

Platonic solids are understood here as the subset of polytopes whose vertices are generated, starting from a single point in $\R^n$, $n\geq3$, by the action of a finite Coxeter group. Platonic polytopes are distinguished by the fact that their faces $f_d$ of any dimension $0\leq d\leq n-1$ are Platonic solids of lower dimension, and that they are transformed into each other by the action of the Coxeter group, i.e., they belong to one orbit of the corresponding Coxeter group.

\begin{figure}[h]
\includegraphics{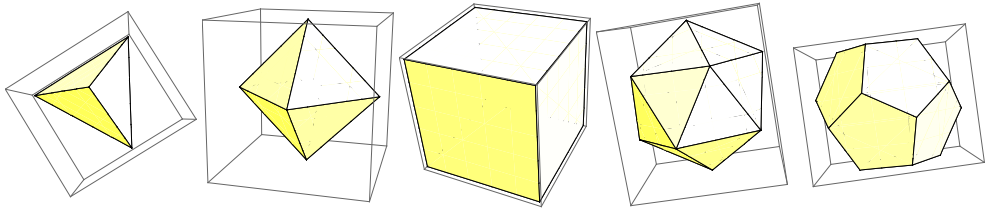}\caption{{\footnotesize The Platonic solids in the 3 dimensional case: tetrahedron, cube, octahedron, icosahedron, dodecahedron.}}\label{plato}
\end{figure}

It has been known since antiquity that there are five platonic solids in $\R^3$~\cite{Cox}, namely, the regular tetrahedron, cube, octahedron, icosahedron and dodecahedron (see Fig.~\ref{plato}). The underlying Coxeter groups\footnote{Reflection groups are denoted  by symbols commonly used for respective simple Lie algebras~\cite{Hum}. Finite Coxeter groups with no connection to Lie algebras are denoted by $H_2$, $H_3$, $H_4$.} are $A_3$, $B_3$, $C_3$ and $H_3$. They are the lowest-dimensional cases we consider by our method.

The root systems of finite Coxeter groups of any type~\cite{D} allow one to extend our method to higher dimensions. For $\R^4$, a classification of Platonic solids was done more than century ago by Schl\"afi~\cite{sh}. In this case, the Coxeter groups are $A_4$, $B_4$, $C_4$, $F_4$ and $H_4$. It is also known that, in any dimension $\geq5$, there are only three such solids generated by groups of types $A_n$, $B_n$ and $C_n$ , namely the simplex, the hypercube and the cross-polytope. They correspond to the regular tetrahedron, cube and octahedron in 3-dimensional space.

For every Platonic solid there exists its dual, which  is also a Platonic solid. We describe both members of each dual pair. For $n\geq3$, the dual pair of $A_n$ consists of two identical solids oriented differently in space. The Platonic solids of $B_n$ and $C_n$ form the dual pair. For $F_4$, $H_3$ and $H_4$, the dual pairs are formed by different solids.

In this paper, we derive these results by new rather simple means, using elementary decoration rules for the corresponding connected Coxeter-Dynkin diagrams~\cite{CKPS}. We also provide a constructive method for building faces $f_d$ of dimensions $0\leq d\leq n-1$. Although we focus on polytopes of dimension $d\geq3$, it is also useful to consider $n=2$ Platonic solids because they occur as 2-dimensional faces of higher-dimensional polytopes.

The general idea of the diagram decoration method \cite{CKPS} is to consider the Coxeter group $W(\mathfrak g)$ of the simple Lie algebra $\mathfrak g$, that is, the symmetry group of a given solid, and to  identify the subgroup $G_{s}(f_d)$ that pointwise stabilizes the given face $f_d$, and the subgroup $G_{f}(f_d)$, which is the symmetry group of the face $f_d$. Then we have
\begin{gather}
G_{s}(f_d)\times G_{f}(f_d)\subset W(\mathfrak g)\,,\qquad 0\leq d\leq n-1\,.
\end{gather}
Since both $G_{s}(f_d)$ and $G_{f}(f_d)$ are generated by reflections defined by some of the simple roots of $\mathfrak g$, it is possible, indeed convenient, to distinguish by different decoration the simple roots of the Coxeter-Dynkin diagram defining the reflections generating the two subgroups. The decorated diagrams allow one to identify all of the polytope faces that belong to different orbits of the corresponding Coxeter groups, and to count how many times each face occurs on the polytope.

Nodes of the Coxeter-Dynkin diagram we decorate by one of the three symbols\footnote{There is a 1-1 correspondence between symbols used to decorate the Coxeter-Dynkin diagrams here and in previous papers. More precisely, $\blacklozenge$, $\square$, $\lozenge$ correspond respectively to $\odot$, $\square$, $\boxtimes$, in \cite{MPvor,MPkal,MPtalk} and to $\bigcirc$, $\Box$, $\boxtimes$ in~ \cite{CKPS}.}, $\blacklozenge$, $\square$, $\lozenge$, according to the rules in Subsection~\ref{rules}.

Decorated Coxeter-Dynkin diagrams are a powerful method of great generality \cite{MPkal}, which could be used to solve other problems, but which so far remain underused. In \cite{MPvor}, the method was used to describe Voronoi and Delone cells in root lattices of all simple Lie algebras $\mathfrak g$. It still has to be used to describe the Voronoi and Delone cells in weight lattices \cite{CS}, and for other problems as well \cite{ful,C70} ...

In this paper, we use a version \cite{CKPS} of the method to describe the Platonic solids together with all of their faces. In dimension 4, equivalent results can be found among the entries of Table~3 of~\cite{CKPS}.

Decoration rules are recursive. Starting from a seed decoration, which provides information about the vertices $f_0$ of the polytope, the procedure consists in modifying the decoration step-by-step. At each step, the modified decoration provides information about the faces of dimension greater by 1. Decoration rules are quite general and are particularly simple when applied to Platonic polytopes. The same set of decorations applies to Coxeter-Dynkin diagrams with the same number of nodes.
As the links between diagram nodes do not affect decoration rules, the links need not be drawn. Only when working with a specific reflection group, the decorated diagram with all links removed is viewed superimposed on the appropriate Coxeter-Dynkin diagram.

In general, the number of copies of the face $f_d$, contained in a given polytope, is equal to the ratio of orders of the Coxeter groups \cite{CKPS},
\begin{equation}\label{facecount}
\# f_d=\frac{|W(\mathfrak g)|}{|W(G_{s}(f_d)|\ |W(G_{f}(f_d)|}\,,
\end{equation}
where $|W|$ is the order of the corresponding reflection group. The subgroups $G_s$ and $G_f$ are read as sub-diagrams of the Coxeter-Dynkin diagram of $\mathfrak g$. Their nodes are identified by the  appropriate decoration for $f_d$, namely $\lozenge$ for reflections generating $G_s$, and $\blacklozenge$ for $G_f$.

The symmetry group of an intersection of two faces, say $f_k$ and $f_j$ of dimensions $k$ and $j$, respectively, is the reflection group $G_f(f_k)\cap G_f(f_j)$ generated by reflections decorated by black rhombi and appearing in the diagram of either face. See the motivating example~\ref{motiv}.

Or, to answer the inverse of this question: Given a face, say $f_k$, of a polytope, how many faces $f_d$ of higher dimension, $k<d$, have $f_k$ in common? Frequently, one would want to know the number of edges originating in the same vertex, $k=0,\ d=1$, or the number of faces $f_d$ meeting in an edge ($k=1$). The answer is the size of the orbit of the stabilizer $\operatorname{Stab}(f_k)$ in $W(\mathfrak g)$ when it is acting on $f_d$. The stabilizer is the reflection group generated by the reflections labeled by rhombi of both colours in the decorated diagram of $f_k$. See examples~\ref{edgmeetver},~\ref{facemeet}.


\section{Preliminaries}

\subsection{Finite reflection groups}\

The well-known results of the classification of finite dimensional simple Lie algebras of any rank and type $n\geq1$ are used here only to identify the finite reflection group $W(\mathfrak g)$, called the Weyl group, or, equivalently, the crystallographic Coxeter group. In addition. we also consider polytopes with symmetries of finite non-crystallographic Coxeter groups denoted by $H_2$, $H_3$, $H_4$, together with their simple root diagrams.

\begin{table}[h]
{\footnotesize
\addtolength{\tabcolsep}{-3pt}
\begin{center}
\begin{tabular}{|c||c|c|c|c|c|c|c|c|c|}
\hline
$\mathfrak g$ &$\ A_n\ $
        &$\ B_n\ $
        &$\ C_n\ $
        &$\ D_n\ $
        &$\ F_4\ $
\\\hline\hline
$|\Delta|$   &$n(n+1)$
             &$2n^2$
             &$2n^2$
             &$2n^2-2n$
             &$48$
\\\hline
$|W|$       &$(n+1)!$
            &$n!\cdot 2^n$
            &$n!\cdot 2^n$
            &$n!\cdot 2^{n-1}$
            &$2^7\cdot 3^2$
\\\hline
\end{tabular}
\bigskip
\begin{tabular}{|c||c|c|c|}
\hline
        &$\ H_2\ $
        &$\ H_3\ $
        &$\ H_4\ $
\\\hline\hline
$|\Delta|$   &$10$
             &$30$
             &$60$
\\\hline
$|W|$       &$2\cdot 5$
            &$2^3 \cdot 3 \cdot 5$
            &$2^6 \cdot 3^2 \cdot 5^2$
\\\hline
\end{tabular}
\bigskip
\caption{{\footnotesize Number $|\Delta|$ of non-zero roots of simple Lie algebras $\mathfrak g$ and the orders $|W|$ of their Weyl group. Number of roots and the orders of the three non-crystallographic Coxeter groups.}} \label{rootnumbers}
\end{center}}
\end{table}

The geometry of the set of simple roots (relative lengths and relative angles) in the real Euclidean space $\R^n$ is described by well-known conventions implied in drawing the corresponding Coxeter-Dynkin diagrams (see, for example,~\cite{Hum}).

Reflections generating the Coxeter groups act in the $n$-dimensional real Euclidean space $\R^n$ spanned by the simple roots $\a_1,\ldots,\a_n$,
according to
\begin{gather}\label{reflection}
r_kx=x-\frac{2\l\a_k,x\r}{\l\a_k,\a_k\r}\a_k\,,\qquad
   x\in\R^n,\quad  k=1,2,\dots,n,
\end{gather}
that is, $r_k$ is the reflection in the hyperplane of dimension $n-1$, containing the origin of $\R^n$, and orthogonal to the simple root $\alpha_k$,
 $k\in\{1,\dots,n \}$.

Instead of $\alpha$-basis of simple roots, we use $\omega$-basis. Two bases are linked by the Cartan matrix $\mathfrak C$, (for example~\cite{Hum, pat})
\begin{equation}
  \alpha = \mathfrak C \omega, \qquad \quad \mathfrak C = \left( \frac{ 2(\alpha_i | \alpha_j )}{(\alpha_j | \alpha_j )}\right), \quad i,j\in \{1,\ldots,n\}.
\end{equation}




The decorations described in this paper should be viewed superimposed on the appropriate Coxeter-Dynkin diagram. Often, the same decoration applies to several diagrams. It should therefore be read in the context of that diagram when some specific information about the face of the polytope needs to be deduced from it.

\subsection{Recursive diagram decoration rules}
\label{rules}\

Decoration rules for connected Coxeter-Dynkin diagrams are recursive. They apply to diagrams of all simple Lie algebras and to any polytope generated by the corresponding Coxeter group starting from a single point in $\R^n$. They can therefore be written without reference to Platonic solids.
\smallskip

\noindent
The `grammar' for every decoration:
\begin{itemize}
\item A node of the diagram can be decorated by $\square$, $\lozenge$, or by $\blacklozenge$.

\item In a diagram of $n$-nodes, there can be up to $n$ squares placed in any node.

\item Any connected pair of nodes must not carry $\lozenge$ and $\blacklozenge$ side by side.
\end{itemize}

\noindent
Assuming that a starting decoration complies with the grammar rules, the decoration describes a face of dimension equal to the number of $\blacklozenge$ in it. For dual polytopes, the same decoration refers to the face of dimensions equal to the number of $\lozenge$.  Usually one starts from vertices, faces of 0-dimension, in which case there is no $\blacklozenge$.
\smallskip

\noindent
Decoration rules:
\begin{enumerate}
\item Replace one of the squares by $\blacklozenge$.

\item Change to $\square$ each $\lozenge$ that became adjacent to the new $\blacklozenge$.

\item Repeat steps (1) and (2) as long as there are any $\lozenge$.
\end{enumerate}


The following proposition is a direct consequence of the decoration rules.

\begin{prop}\

A polytope has one Coxeter group orbit of faces $f_d$ for every dimension $0\leq d\leq n-1$, i.e. it is Platonic, provided
\newline
(i)\ The Coxeter-Dynkin diagram forms a connected line with no branches or loops.
\newline
(ii)\ The seed decoration has one square placed at either of the extreme nodes of the Coxeter-Dynkin diagram.
\end{prop}

\begin{coro}\label{coxdiag}\

Coxeter-Dynkin diagrams that give rise to Platonic solids in $\R^n$, where $n\geq 3$, are of types

{\footnotesize
\parbox{.6\linewidth}{\setlength{\unitlength}{2pt}
\def\kr{\circle{4}}
\def\cr{\circle*{4}}
\thicklines
\begin{picture}(140,30)
\put( 0,22){\makebox(0,0){${A_n}$}}
\put(10,22){\kr}
\put(20,22){\kr}
\put(40,22){\kr}
\put(27,21.7){$\dots$}
\put(50,22){\kr}
\put(12,22){\line(1,0){6}}
\put(42,22){\line(1,0){6}}
\put(22,22){\line(1,0){4}}
\put(34,22){\line(1,0){4}}
\put( 0,12){\makebox(0,0){${B_n}$}}
\put(10,12){\kr}
\put(20,12){\kr}
\put(27,11.7){$\dots$}
\put(40,12){\kr}
\put(50,12){\cr}
\put(12,12){\line(1,0){6}}
\put(22,12){\line(1,0){4}}
\put(34,12){\line(1,0){4}}
\put(42,13){\line(1,0){6}}
\put(42,11){\line(1,0){6}}
\put( 0,2){\makebox(0,0){${C_n}$}}
\put(10,2){\cr}
\put(20,2){\cr}
\put(40,2){\cr}
\put(50,2){\kr}
\put(12,2){\line(1,0){6}}
\put(22,2){\line(1,0){4}}
\put(27,1.7){$\dots$}
\put(34,2){\line(1,0){4}}
\put(42,3){\line(1,0){6}}
\put(42,1){\line(1,0){6}}
\put(102,22){\makebox(0,0){${F_4}$}}
\put(110,22){\kr}
\put(120,22){\kr}
\put(130,22){\cr}
\put(140,22){\cr}
\put(112,22){\line(1,0){6}}
\put(122,23){\line(1,0){6}}
\put(122,21){\line(1,0){6}}
\put(132,22){\line(1,0){6}}
\put(102,12){\makebox(0,0){${H_3}$}}
\put(110,12){\kr}
\put(120,12){\kr}
\put(130,12){\kr}
\put(124,13){$5$}
\put(112,12){\line(1,0){6}}
\put(122,12){\line(1,0){6}}
\put(102,2){\makebox(0,0){${H_4}$}}
\put(110,2){\kr}
\put(120,2){\kr}
\put(130,2){\kr}
\put(140,2){\kr}
\put(134,3){$5$}
\put(112,2){\line(1,0){6}}
\put(122,2){\line(1,0){6}}
\put(132,2){\line(1,0){6}}
\end{picture}}}
\end{coro}

\smallskip

\subsection{Dual polytopes.}\

In this paper, we define dual polytopes and their faces using decorated Coxeter-Dynkin diagrams:
\newline
The role of decoration elements, $\lozenge$ and $\blacklozenge$, are reversed in the dual polytope. That is, $\lozenge$ specifies reflections that generate the symmetry group $G_f$ of the dual face, while $\blacklozenge$ provides the generating reflections of the subgroup $G_s$ that stabilizes the face pointwise.

The dual polytope of a Platonic polytope is also platonic. Specifically, we have dual pairs of platonic polytopes in $\R^n$, $n\geq3$,

\begin{align}
A_n: \quad&
\textrm{Both polytopes coincide, but are differently oriented in } \R^n \notag\\
B_n: \quad&
B_n \textrm{ and } C_n \textrm{ polytopes are dual to each other in } \R^n \notag\\
F_4: \quad&
\textrm{The dual polytopes in } \R^4 \textrm{ are different (see Table~\ref{edges4})} \notag \\
H_3: \quad&
\textrm{The dual polytopes are the icosahedron and dodecahedron in } \R^3 \notag\\
H_4: \quad &
\textrm{The dual polytopes are different in } \R^4 \textrm{ (see Table~\ref{edges4})} \notag
\end{align}

A particular decoration that carries information about a face $f_d$ of a polytope in $\R^n$ can be read as  information about the face of the dual polytope, say $\tilde f_{n-d-1}$.

\section{Platonic solids}

\subsection{Motivating example: Platonic solids in dimension 3}\label{motiv}\

Let us illustrate the decoration method on the transparent example before presenting the general rules. Consider the classical platonic solids in $\R^3$. Connected Coxeter-Dynkin diagrams of the groups, generated by three reflections, are of types $A_3$, $B_3$, $C_3$ and $H_3$ (see Corollary~\ref{coxdiag}).

We assume that diagram nodes are numbered 1,2,3 left to right, that each node represents the corresponding reflection $r_1$, $r_2$, $r_3$ acting in $\R^3$. Decorations described in this subsection should be viewed superimposed on any of the four diagrams (Cor.~\ref{coxdiag}).

Let the seed point (face $f_0$) be either $\omega_1$ or $\omega_3$, corresponding to the dual pair of polytopes. We denote by $\square$ the reflections that move the seed point and by $\lozenge$ the reflections that stabilize it. Consider two seed point decorations of all four diagrams (Cor.~\ref{coxdiag}), namely
\begin{gather}
\square\ \lozenge\ \lozenge\qquad\text{seed point $\omega_1$}
\label{first} \notag \\
\lozenge\ \lozenge\ \square\qquad\text{seed point $\omega_3$} \notag
\label{second}
\end{gather}
Here, $\square$ stands for the seed point reflection. The sub-diagram decorated by $\lozenge\ \lozenge$ identifies reflections generating $G_{s}$. The symmetry group $G_{f}(f_0)=1$ of the face $f_0$ is trivial, the face is just a point.

The orders of reflection groups $A_3, B_3 ,H_3$ in $\R^3$ are shown in Table~\ref{rootnumbers}.

The groups $W(G_{s})$ stabilizing $\omega_1$ are respectively
\begin{gather}
|W(A_2)|=6,\qquad
|W(B_2)|=|W(C_2)|=8,\qquad
|W(H_2)|=10
\end{gather}

Based on the number of vertices, we have a tetrahedron for $A_3$, an octahedron for $B_3$ and $C_3$, and a dodecahedron for $H_3$.

The groups $W(G_{s})$ stabilizing $\omega_3$ are all of type $A_2$. Hence, the numbers of vertices are those of a tetrahedron for $A_3$, a cube for $B_3$ and $C_3$, and a icosahedron for $H_3$. The dual tetrahedra of $A_3$ are differently oriented in $\R^3$, but otherwise are identical.

The second step describes the edges of the solids. The decorations are
\begin{gather}
\blacklozenge\ \square\ \lozenge \notag\\
\lozenge\ \square\ \blacklozenge\ \notag
\end{gather}

The symmetry group generated by $r_1$ refers to the edge with end points $\omega_1$ and $r_1\omega_1$. For the dual polytope, the symmetry group is generated by $r_3$, and it refers to the edge with end points $\omega_3$ and $r_3\omega_3$. Note that this description of edges is independent of the Coxeter groups used in~Cor.~\ref{coxdiag}.

In all cases, the group $G_s(f_1)\times G_f(f_1)$ is $W(A_1\times A_1)$ of order 4. Consequently, the number of edges is $|W(\mathfrak g)|/4$. We get 6 for $A_3$ (tetrahedron), 12 for $B_3$ and $C_3$ (cube and octahedron) and 30 for $H_3$ (icosahedron and dodecahedron).

The final step in the decoration
\begin{gather}
\blacklozenge\ \blacklozenge\ \square \label{f2o1}\\
\square\ \blacklozenge\ \blacklozenge\ \label{f2o2}
\end{gather}
describes the 2-dimensional faces $f_2$. The sub-diagram $\blacklozenge\ \blacklozenge$ in \eqref{f2o1} is the group $G_f(f_2)=W(A_2)$ of order $|W(A_2)|=6$ for all cases.
$$
\#f_2\ \ :\quad\frac{|W(A_3)|}{|W(A_2)|}=4\,,\qquad
                \frac{|W(B_3)|}{|W(A_2)|}=
                \frac{|W(C_3)|}{|W(A_2)|}=8\,,\qquad
                \frac{|W(H_3)|}{|W(A_2|)}=20\,.
$$
The same sub-diagram in \eqref{f2o2} stands for $G_f(f_2)=W(A_2)$ in the case of $A_3$, $G_f(f_2)=W(C_2)$ for $B_3$ and $C_3$, and $G_f(f_2)=W(H_2)$ for $H_3$. Therefore
$$
\#\tilde{f}_2\ \ :\quad\frac{|W(A_3)|}{|W(A_2)|}=4\,,\qquad
                \frac{|W(B_3)|}{|W(C_2)|}=
                \frac{|W(C_3)|}{|W(C_2)|}=6\,,\qquad
                \frac{|W(H_3)|}{|W(H_2|)}=12\,.
$$
The shape of face $f_2$ can be easily determined from the relative angles of the mirrors $r_1$ and $r_2$ acting on $\omega_1$ or from reflections $r_2$ and $r_3$ acting on $\omega_3$.

\smallskip

A summary of properties of Platonic solids in $\R^3$ and their faces is shown in Table~\ref{platoDim3}.

\begin{table}[h]
{\footnotesize
$  \begin{array}{cccccccccc}
\textrm{Face} & G_f & G_s &d & v&N(A_3)&N(B_3)&N(H_3) & \textrm{Platonics}\\ \hline
\square\ \lozenge\ \lozenge  & 1&r_2,r_3&0&2&4&6&12&\begin{array}{cc}
                                     \surd&
                                   \end{array}\\
\lozenge\ \lozenge\ \square & 1&r_1,r_2&0&2&4&8&20&\begin{array}{cc}
                                     &\surd
                                   \end{array}\\
\blacklozenge\ \square\ \lozenge & r_1&r_3&1&1&6&12&30&\begin{array}{cc}
                                     \surd&
                                   \end{array}\\
\lozenge\ \square\ \blacklozenge & r_3&r_1&1&1&6&12&30&\begin{array}{cc}
                                     &\surd
                                   \end{array}\\
\blacklozenge\ \blacklozenge\ \square & r_1,r_2&1&2&0&4&8&20&\begin{array}{cc}
                                     \surd&
                                   \end{array}\\
\square\ \blacklozenge\ \blacklozenge\ & r_2,r_3&1&2&0&4&6&12&\begin{array}{cc}
                                     &\surd
                                   \end{array}\\
  \end{array}
$
\bigskip
\caption{{\footnotesize Each line describes one face $f_d$ in the Platonic solids, with either of the symmetry groups $A_3$, $B_3$, $H_3$. Reflections generating $G_f$ and $G_s$ are in the second and third column. Column $d$ contains the dimension of the face, column $v=2-d$ shows the dimension of the corresponding face of the dual Platonic polytope. Columns $N(A_3)$, $N(B_3)$ and $N(H_3)$ contain the number of times the face $f_d$ occurs in the polytope. Vertically aligned marks in the last column indicate which faces belong to the same polytope of the dual pair.}} \label{platoDim3}}
\end{table}


\begin{example}\

The dual pair of platonic solids of $A_3$ consists of two tetrahedra differently oriented in $\R^3$, their seed points being $\omega_1$ and $\omega_3$.

The vertices of the two tetrahedra are the following
\begin{gather}
\omega_1,\ -\omega_1+\omega_2,\ -\omega_2+\omega_3,\ -\omega_3 \notag\\
\omega_3,\ -\omega_3+\omega_2,\ -\omega_2+\omega_1,\ -\omega_1 \notag
\end{gather}

\end{example}


\subsection{Platonic solids in dimension 4}\label{solid4dim}\

Connected Coxeter-Dynkin diagrams of the groups generated by four reflections are of types $A_4$, $B_4$, $C_4$, $F_4$ and $H_4$ (see (Corollary~\ref{coxdiag})).

Using table~\ref{platoDim3} conventions, properties of the 4-dimensional platonic solids, together with their duals, are summarized in Table~\ref{platoDim4}.

\begin{table}[h]
{\footnotesize
$  \begin{array}{ccccccccccc}
\textrm{Face} & G_f & G_s &d & v&N(A_4)&N(B_4)&N(F_4)&N(H_4) & \textrm{Platonics}\\ \hline
\square\ \lozenge\ \lozenge\ \lozenge & 1&r_2,r_3,r_4&0&3&5&8&24&120&\begin{array}{cc}
                                     \surd&
                                   \end{array}
\\
\lozenge\ \lozenge\ \lozenge\ \square & 1&r_1,r_2,r_3&0&3&5&16&24&600&\begin{array}{lr}
                                     &\surd
                                   \end{array}\\
\blacklozenge\ \square\ \lozenge\ \lozenge & r_1&r_3,r_4&1&2&10&24&96&720&\begin{array}{cc}
                                     \surd&
                                   \end{array}\\
\lozenge\ \lozenge\ \square\ \blacklozenge & r_4&r_1,r_2&1&2&10&32&96&1200&\begin{array}{lr}
                                     &\surd
                                   \end{array}\\
\blacklozenge\ \blacklozenge\ \square\ \lozenge & r_1,r_2&r_4&2&1&10&32&96&1200&\begin{array}{lr}
                                     \surd&
                                   \end{array}\\
\lozenge\ \square\ \blacklozenge\ \blacklozenge & r_3,r_4&r_1&2&1&10&24&96&720&\begin{array}{lr}
                                     &\surd
                                   \end{array}\\
\blacklozenge\ \blacklozenge\ \blacklozenge\ \square & r_1,r_2,r_3&1&3&0&5&16&24&600&\begin{array}{lr}
                                     \surd&
                                   \end{array}\\
\square\ \blacklozenge\ \blacklozenge\ \blacklozenge & r_2,r_3,r_4&1&3&0&5&8&24&120&\begin{array}{lr}
                                     &\surd
                                   \end{array}\\
  \end{array}
$
\bigskip
\caption{{\footnotesize Underlying reflection groups in $\R^4$ are of types $A_4$, $B_4$, $F_4$, $H_4$. For conventions, see the caption of Table~\ref{platoDim3}.
}} \label{platoDim4}}
\end{table}

\begin{example}\label{edgmeetver}\

Let us answer the question posed in the Introduction: How many edges $f_1$ meet in a vertex $f_0$ of any Platonic solid?

Decorations of Coxeter-Dynkin diagrams for 4-dimensional Platonic solids and their duals are:

\begin{gather}
f_0: \quad  \square\ \lozenge\ \lozenge\ \lozenge  \qquad\qquad \qquad f_1: \quad  \blacklozenge\ \square\ \lozenge\ \lozenge \label{f1}\\
\tilde{f}_0: \quad  \lozenge\ \lozenge\ \lozenge\ \square \qquad\qquad \qquad \tilde{f}_1: \quad  \lozenge\ \lozenge\  \square\ \blacklozenge \label{f'1}
\end{gather}

From table~\ref{platoDim4}, we read the stabilizers of faces

\begin{gather}
G_s({f_0}) = \l r_2,r_3,r_4 \r \qquad \qquad  G_s(f_1) = \l r_3,r_4 \r \label{Gf1} \\
G_s(\tilde{f}_0) = \l r_1,r_2,r_3 \r \qquad \qquad G_s(\tilde{f}_1) = \l r_1,r_2 \r \label{Gf'1}
\end{gather}

The edges originating in $f_0$ or $\tilde{f}_0$ are generated by the stabilizer of $f_0$ or $\tilde{f}_0$, respectively. 

They are equal to $G_s(f_1)$ for~\eqref{Gf1}, and $G_s(\tilde{f}_1)$ for~\eqref{Gf'1}.

The formula for the number of edges meeting in one vertex can be written as
 \begin{equation}
   \# f_1(f_0) = \frac{|W(G_s(f_0))|}{|W(G_s(f_1))|},
 \end{equation}

or more generally, if we consider $f_{d}(f_{d-1})$ of faces $f_d$ having in common faces $f_{d-1}$ for $1\le d \le n-2$:

\begin{prop}\label{numvert}\

The number of faces $f_d$ meeting in a face $f_{d-1}$ for $d\in \{1,...n-2\}$ is equal to the ratio of orders of Coxeter groups of the subgroups $G_s(f_d), G_s(f_{d-1})$ which stabilize a given face $f_d, f_{d-1}$, respectively
 \begin{equation}
   \# f_{d}(f_{d-1}) = \frac{|W(G_s(f_{d-1}))|}{|W(G_s(f_{d}))|}.
 \end{equation}
\end{prop}

\bigskip

The groups $W(G_s(f_1))$ for~\eqref{f1} are respectively $W(A_2), W(B_2), W(A_2)$ and $W(H_2)$, and for~\eqref{f'1} $W(G_s(\tilde{f}_1))=W(A_2)$.

Table~\ref{edges4} summarizes this proposition.

\begin{table}[h]
{\centering \footnotesize
$\begin{array}{cccc}
  \textrm{Name of Polytope} & \# f_0 & \# f_1(f_0) & \# f_2(f_1) \\ \hline
  \textrm{Pentatope} & 5 & 4 & 3\\
  \textrm{16-cell} & 8 & 6 & 4\\
  \textrm{Tessaract} & 16 & 4 & 3\\
  \textrm{24-cell} & 24 & 8 & 3\\
  \textrm{600-cell} & 120 & 20 & 5\\
  \textrm{120-cell} & 600 & 4 & 3
\end{array}$
\bigskip\caption{{\footnotesize The first column contains the names of the 4-dimensional Platonic polytopes and their duals. The second column shows the number of vertices in the polytope. The underlying symmetry group for each line of the table can be identified by comparing the number of vertices $\# f_0$ with the corresponding entries in Table~\ref{platoDim4}. The third column contains the number $\# f_1(f_0)$ of edges meeting at each vertex, and the last column contains the number $\# f_2(f_1)$ of 2-dimensional faces meeting at each edge.}}\label{edges4}}
\end{table}

Both the pentatope and 24-cell are self-dual, the 16-cell is dual to the tessaract, and the 600-cell and 120-cell are dual to each other~\cite{Cox}.

\end{example}

\bigskip

\subsection{Platonic solids in dimension $\ge 5$}\

\smallskip
The are only 3 types of platonic solids for dimension more than 4, namely the simplex, hypercube and cross-polytope. The simplex is self-dual, and the cross-polytope and hypercube are dual to each other. The corresponding Coxeter-Dynkin diagrams (Cor.~\ref{coxdiag}) are of types $A_n, B_n$ and $C_n$.

\smallskip

Using Table~\ref{platoDim3} conventions, the properties of Platonic solids of dimension $\geq5$, together with their duals, are summarized in Table~\ref{platoDimn}.

\begin{table}[h]
{\footnotesize
$  \begin{array}{ccccccccc}
\textrm{Face} & G_f & G_s &d & v&N(A_n)&N(B_n) & \textrm{Platonics}\\ \hline
\square\ \lozenge\ \lozenge\ \ldots\ \lozenge\ \lozenge & 1&r_2,\ldots,r_n&0&n-1&\tfrac{(n+1)!}{n!}&\tfrac{2^nn!}{2^{n-1}(n-1)!}&\begin{array}{cc}
                                     \surd&
                                   \end{array}
\\
\lozenge\ \lozenge\ \lozenge\ \ldots\ \lozenge\ \square & 1&r_1,\ldots,r_{n-1}&0&n-1&\tfrac{(n+1)!}{n!}&\tfrac{2^nn!}{n!} &\begin{array}{lr}
                                     &\surd
                                   \end{array}\\
\blacklozenge\ \square\ \lozenge\ \ldots\ \lozenge\ \lozenge & r_1&r_3,\ldots,r_n&1&n-2&\tfrac{(n+1)!}{2!(n-1)!}&\tfrac{2^nn!}{2!2^{n-2}(n-2)!}&\begin{array}{cc}
                                     \surd&
                                   \end{array}\\
\lozenge\ \lozenge\ \lozenge\ \ldots\ \square\ \blacklozenge & r_n&r_1,\ldots,r_{n-2}&1&n-2&\tfrac{(n+1)!}{(n-1)!2!}&\tfrac{2^nn!}{(n-1)!2!}&\begin{array}{lr}
                                     &\surd
                                   \end{array}\\
\blacklozenge\ \blacklozenge\ \square\ \ldots\ \lozenge\ \lozenge & r_1,r_2&r_4,\ldots,r_n&2&n-3&\tfrac{(n+1)!}{3!(n-2)!}&\tfrac{2^nn!}{3!2^{n-3}(n-3)!}&\begin{array}{lr}
                                     \surd&
                                   \end{array}\\
\lozenge\ \lozenge\ \square\ \ldots\ \blacklozenge\ \blacklozenge & r_{n-1},r_n&r_1,\ldots,r_{n-3}&2&n-3&\tfrac{(n+1)!}{(n-2)!3!}&\tfrac{2^nn!}{(n-2)!2^2 2!}&\begin{array}{lr}
                                     &\surd
                                   \end{array}\\
                             \vdots    &\vdots&\vdots&\vdots   &\vdots&\vdots&\vdots&\vdots\\
\blacklozenge\ \blacklozenge\ \blacklozenge\ \ldots\ \square\ \lozenge & r_1,\ldots,r_{n-2}&r_n&n-2&1&\tfrac{(n+1)!}{(n-1)!2!}&\tfrac{2^nn!}{2!(n-1)!}&\begin{array}{lr}
                                     \surd&
                                   \end{array}\\
\lozenge\ \square\ \blacklozenge\ \ldots\ \blacklozenge\ \blacklozenge & r_3,\ldots,r_n&r_1&n-2&1&\tfrac{(n+1)!}{2!(n-1)!}&\tfrac{2^nn!}{2!2^{n-2}(n-2)!}&\begin{array}{lr}
                                     &\surd
                                   \end{array}\\
\blacklozenge\ \blacklozenge\ \blacklozenge\ \ldots\ \blacklozenge\ \square & r_1,\dots,r_{n-1}&1&n-1&0&\tfrac{(n+1)!}{n!}&\tfrac{2^nn!}{n!}&\begin{array}{lr}
                                     \surd&
                                   \end{array}\\
\square\ \blacklozenge\ \blacklozenge\ \ldots\ \blacklozenge\ \blacklozenge & r_2,\ldots,r_n&1&n-1&0&\tfrac{(n+1)!}{n!}&\tfrac{2^nn!}{2^{n-1}(n-1)!}&\begin{array}{lr}
                                     &\surd
                                   \end{array}\\
  \end{array}
$
\bigskip
\caption{{\footnotesize Underlying reflection groups in $\R^n$, $(n\geq5)$ are of types $A_n$, $B_n$. 
}} \label{platoDimn}}
\end{table}

Table~\ref{meetingfaces} summarizes this proposition~\ref{numvert} for any n-dimensional polytope.

\begin{table}[h]
{\centering \footnotesize
$\begin{array}{ccccccc}
  \textrm{Name of Polytope} & \# f_0 & \# f_1(f_0) & \# f_2(f_1) & \cdots & \# f_{n-3}(f_{n-4}) &\# f_{n-2}(f_{n-3}) \\ \hline
  \textrm{Simplex ($(n+1)$-cell)} & n+1 & n & n-1 & \cdots & 4 & 3 \\
  \textrm{Cross-polytope ($2^n$-cell)} & 2n & 2(n-1) & 2(n-2) & \cdots & 8 & 4 \\
  \textrm{Hypercube ($2n$-cell)} & 2^n & n & n-1 & \cdots & 4 & 3
\end{array}$
\bigskip\caption{{\footnotesize The first column contains the names of the n-dimensional Platonic polytopes and their duals. The second column shows the number vertices in the polytope, the remaining columns contain the number of faces $f_{d}$ meeting at each face $f_{d-1}$ for $d\in \{1,\ldots,n-2\}$, $n\ge5$.}}\label{meetingfaces}}
\end{table}

Analog to proposition~\ref{numvert}, we can give the formula for finding the number of faces $f_d$ meeting at face $f_c$, where $c<d$.

\begin{prop}\

The number of faces $f_d$ meeting at face $f_{c}$ for $0\le c<d<n-2$ is equal the ratio of orders of Weyl groups of subgroups $G_s(f_c), G_s(f_{d})$, which stabilize a given face $f_c$ and $f_{d}$, respectively
 \begin{equation}
   \# f_{d}(f_{c})
    = \frac{|W(G_s(f_{c}))|}{|W(G_s(f_{d}))|}.
 \end{equation}
\end{prop}
\medskip

\begin{example}\label{facemeet}\

We now consider the 4-dimensional faces appearing in 5-dimensional polytopes.

Diagrams of the face $f_4$ are:

$$\blacklozenge\ \blacklozenge\ \blacklozenge\ \blacklozenge\ \square $$
$$\square\ \blacklozenge\ \blacklozenge\ \blacklozenge\ \blacklozenge$$

The shape of face $f_4$ is determined from the relative angles of mirrors $r_1, \ldots, r_4$ acting on $\omega_1$, and from reflections $r_2, \ldots, r_5$ acting on $\omega_5$. Using Tables~\ref{platoDim4} and \ref{edges4}, one can see from sub-diagram $\blacklozenge\ \blacklozenge\ \blacklozenge\ \blacklozenge$ that, for $A_5$, faces are pentatopes, and for $B_5$, either 16-cells or tessaracts.
\end{example}

\section{Concluding remarks}

A seed decoration is often applied to several Coxeter diagrams. See an example in~\ref{solid4dim}. The same is true for recursive decorations that follow from the seed. Therefore, the set of such decorations could be viewed as a `generic polytope' pertinent to all  Coxeter groups with the same connectivity of their Coxeter-Dynkin diagrams. Can anything be learned from the generic set of decorations?

In this paper, we focused on Coxeter groups whose diagrams are connected. In general, the decoration rules can be used for polytopes of Lie algebras that are semisimple but not simple, i.e. their Coxeter-Dynkin diagrams are disconnected. In this case, each connected component must have its own seed point, indicated in the initial decoration. Consequently, there are several orbits of edges, hence, such polytopes are never Platonic.

In \cite{CKPS}, the authors considered semiregular polytopes in 3 and 4 dimensions. It would be interesting to describe semiregular polytopes in a higher dimension.

In recent years evidence has been obtained that there exist in nature molecules with imperfect symmetries \cite{ful, C70} of $W(\mathfrak g)$ type, not necessarily Platonic. They could be considered as $W(\mathfrak g)$ symmetries broken to a subgroup.

\subsection*{Acknowledgements}\

\noindent
The author would like to express her gratitude to the Centre de recherches math\'ematiques, Universit\'e de Montr\'eal, for the hospitality extended to her during her postdoctoral fellowship.

\noindent
She would like also to thank Dr. J. Patera for stimulating discussions and comments.

\noindent
She is grateful to the MIND Reseach institute, Santa Ana, Calif., to MITACS and to OODA Technologies for partial support.



\end{document}